\documentclass[12pt,a4paper]{article}

\usepackage{amsfonts,amssymb,amsmath,latexsym}
\usepackage{graphicx}
\usepackage{vmargin}

\title {On the number of non-isomorphic subspaces of a Banach
space.}
\author {Valentin Ferenczi and Christian Rosendal}
\date {December 2001}
\newcommand {\ca} {2^\omega}

\newcommand {\N}{\mathbb N}
\newcommand {\Q}{\mathbb Q}
\newcommand {\R}{\mathbb R}

\newcommand{\om}{\omega}

\newcommand{\eps}{\epsilon}

\newcommand {\iso} {\cong}

\newcommand {\tom} {\emptyset}

\newcommand{\norm}[1]{\|{#1}\|}

\newcommand {\Del}{ \; \Big| \;}
\newcommand {\del}{ \; \big| \;}

\newcommand {\mgdv}{\big\{}

\newcommand {\mgdh}{\big\}}

\newcommand {\Intv}{\Big[}
\newcommand {\intv}{\big[}

\newcommand {\Inth}{\Big]}
\newcommand {\inth}{\big]}

\newcommand {\for}{\bigcup}

\newcommand {\snit}{\bigcap}

\newcommand {\og}{\; \land \;}

\newcommand {\go} {\mathfrak}
\newcommand {\ku} {\mathcal}

\newcommand {\e} {\exists}
\renewcommand {\a} {\forall}

\renewcommand {\pmb} {\boldsymbol}
\newcommand{\fed}{\boldsymbol}

\newcommand{\pf}{

\smallskip

\noindent {\it Proof : }}

\newcommand{\pff}{$\hfill  \Box$

\smallskip }

\newtheorem{thm}{Theorem}
\newtheorem{cor}[thm]{Corollary}
\newtheorem{lemme}[thm]{Lemma}
\newtheorem{prop} [thm] {Proposition}
\newtheorem{defi} [thm] {Definition}

\newtheorem{ex}[thm]{Example}

\begin{document}
\maketitle

\begin{abstract}{We study the number of non-isomorphic subspaces of a given
Banach space. Our
main result is the following: let $\go X$ be a Banach space with an
unconditional
basis $\{e_i\}_{i \in \N}$; then either there exists a perfect set $\ku P$
of infinite subsets of $\N$ such that for any two distinct $A,B \in \ku
P$, $[e_i]_{i \in A} \ncong [e_i]_{i \in B}$, or for a residual set
of infinite subsets $A$ of $\N$, $[e_i]_{i \in A}$ is isomorphic to
$\go X$, and in that case,
$\go X$ is isomorphic to its square, to its hyperplanes,
uniformly isomorphic to $\go X \oplus [e_i]_{i \in D}$
for any
 $D\subset \N$ and to a denumerable Schauder decomposition into
 uniformly isomorphic copies of itself.}
\end{abstract}

The starting point of this article is the so called
 "homogeneous space problem", due to S. Banach. It
has been solved at the end of last century by W.T. Gowers (\cite{G1},\cite{G2}),
 R. Komorowski- N. Tomczak (\cite{KT1},\cite{KT2}) and  W.T. Gowers- B.
Maurey (\cite{GM});
 see e.g. \cite{T} for a survey.
Recall that a Banach space is said to be {\em homogeneous} if it is isomorphic
 to all its (infinite-dimensional closed) subspaces. The previously
named authors showed that $l_2$
is the only homogeneous Banach space.
A very natural question was posed to us by G. Godefroy: if a Banach
space is
not isomorphic to $l_2$, then how many mutually non-isomorphic subspaces
must it
contain (obviously, at least $2$).
In this article, we concentrate on spaces with a basis and subspaces of it
spanned by subsequences. We shall also be interested in the relation
of equivalence of basic sequences. By "many" we shall mean the Cantor
concept of cardinality, and sometimes finer concepts from the theory
of classification of equivalence relations in descriptive set theory.

\section{Basic notions about basic sequences.}

Let $\go X$ be some separable Banach space and $\mgdv e_i \mgdh_\om$ a
non zero  sequence in $\go X$. We say that  $\mgdv e_i \mgdh_\om$ is a
basis for $\go X$ if any vector $x$ in $\go X$ can be uniquely written as a
norm convergent series $x=\sum_\omega a_ie_i$.
In that case the functionals $e^*_k(\sum_\omega a_ie_i):=a_k$ are in fact
continuous, as are the projections $P_n(\sum_\omega a_ie_i):=\sum_{i=0}^n
a_ie_i$ and furthermore their norms are uniformly bounded
$\sup_\om\|P_n\|<\infty$; the supremum is called the basic constant of
$\mgdv e_i \mgdh_\om$ and is denoted by $bc(\mgdv e_i \mgdh_\om)$.
If $\mgdv e_i \mgdh_\om$ is some non zero sequence that is a basis for its
closed linear span, written $\intv e_i \inth_\om$, we say that it is a
basic sequence in $\go X$ and its basic constant is defined as before.
The property of $\mgdv e_i \mgdh_\om$ being a basic sequence can also
equivalently be stated as the existence of a constant $K\geq 1$ such that
for any $n\leq m$ and $a_0,a_1,\ldots,a_m\in \R$:
$$\|\sum_{i=0}^n a_ie_i\|\leq K\|\sum_{i=0}^m a_ie_i\|$$
The infimum of such $K$ will then be the basic constant $bc(\mgdv e_i
\mgdh_\om)$.
Suppose furthermore that for any $x=\sum_\omega a_ie_i$ the series actually
converges unconditionally, ie. for any permutation $\sigma$ of $\om$ the
series $\sum_\omega a_\sigma(i)e_\sigma(i)$ converge to $x$, then the basic
sequence is said to be unconditional.
We note that a classical result of Dvoretzsky and Rogers states that in an
infinite dimensional Banach space, in contradiction to $\R$, there is
always an unconditionally convergent series that does not converge
absolutely (ie. $\sum_\om\|x_i\|=\infty$).
Again being an unconditional basis for some closed subspace (which will be
designated by `unconditional basic sequence') is equivalent to there
being a constant $K\geq 1$, such that for all n and $A\subset \mgdv
0,\ldots,n\mgdh, \; a_0,\ldots,a_n\in \R$
$$\||\sum_{i\in A} a_ie_i\||\leq K\|\sum_{i=0}^n a_ie_i\|$$
We will in general only work with normalized basic sequences, ie.
$\|e_i\|\equiv 1$, which always can be obtained by taking
$e'_i:=\frac{e_i}{\|e_i\|}$.
Moreover it will often be useful to to suppose the basis monotone, that is
that the above inequalities holds for $K=1$. This can also be got by
renorming the space by
$$\|\!| \sum_\om a_ie_i\|\!|:=\sup_\om\|\sum_{i=0}^n a_ie_i\|\leq
bc\{e_i\}\|\sum_\om a_ie_i\|$$
in the case of a (maybe) conditional basis and
$$\|\!|\sum_\om a_ie_i\|\!|:=\sup_{<\eps_i>\in \{-1,1\}^\om}\|\sum_\om
\eps_ia_ie_i\|$$
for an unconditional basis.

Two  sequences $(e_i)$ and $(f_i)$
are said to be equivalent if there exists
 $C$ and $a,b$ with $ab \leq C$ such that
$(1/a)\norm{\sum \lambda_i e_i} \leq
\norm{\sum \lambda_i f_i} \leq b\norm{\sum \lambda_i e_i}$, in which case they
are said to $C$-equivalent. For basic sequences, this is equivalent to
saying that for any
choice of reals $(\lambda_i)_{i \in \N}$,
$\sum \lambda_i e_i$ converges if and only if $\sum \lambda_i f_i$ converges.
We shall write $(e_i) \approx (f_i)$ to mean that $(e_i)$ and $(f_i)$ are
equivalent.

A space $\go Y$ is
isomorphic to $\go Z$ if there is a bijective continuous linear map $T$ from
$\go Y$ onto $\go Z$ such
that $T^{-1}$ is also continuous.
They are $C$-isomorphic if this happens for some $T$ such that
$\norm{T}\norm{T^{-1}} \leq C$.
The Banach-Mazur distance between $\go Y$ and $\go Z$, denoted by
$d_{BM}(\go Y,\go Z)$ is defined as the infimum of the $C$'s such that
$\go Y$ and $\go Z$ are $C$-isomorphic.
We shall write ${\go Z} \simeq {\go Y}$ to mean that
$\go Y$ and $\go Z$ are isomorphic,
$\go Z \simeq_K \go Y$ to mean that they are $K$-isomorphic.

\section {Models of separable Banach spaces and their basic sequences.}

To study separable Banach spaces by topological means we need in some way
to make a space out of them. So we turn to descriptive set theory for the
basic tools.
Let $\go X$ be a Polish space and $F(\go X)$ denote the set of closed
subsets of $\go X$. We endow $F(\go X)$ with the following $\sigma$-algebra
that renders it a standard Borel space.
The generators are the following sets, where $U$ varies over the open
subsets of $\go X$:
$$\mgdv F \in F ( \go X) \del F \cap U \neq \tom \mgdh$$
The resulting measurable space is called the \emph{Effros Borel space} of
 $F(\go X)$.

A theorem due to Kuratowski/Ryll-Nardzewski states that there is a sequence
of Borel functions $d_n:F(\go X) \longrightarrow \go X$, such that for
nonempty $F\in F(\go X)$ the set $\mgdv d_n(F) \mgdh $ is dense in $F$.
Supposing now that $\go X $ is a separable Banach space,  we can in a Borel
manner express that $F\in F(\go X)$ isa linear subspace of $\go X$:
$$0\in F \land \a n,m \; \a p,q \in \Q (pd_n(F)+qd_m(F)\in F)$$
Here we have in all honesty to admit that we use the fact that the
relations $F_0 \subset F_1$ and $x\in F$ are Borel in $F(\go X)^2$ and $\go
X \times F(\go X)$ respectively.

But then it is possible to construct a standard Borel space consisting of
all separable (real if you wish) Banach spaces, simply take any
isometrically universal separable Banach space (to be  definite take
$C(\ca)$) and let $\go A$ be the set of closed linear subspaces of it. Call
$\go A$ the Borel space of separable Banach spaces.

There is now a space in which we can express the relations of (linear)
isomorphism,
 isometry, etc. and we will see that their descriptive complexity are as they
should be, eg. analytic for isomorphism.

If we restrict ourselves to certain types of subspaces, the situation is
far less involved. Spaces spanned by subsequences of a
given basis can be identified with
$\ca$.
Also, following Gowers in \cite{G1},\cite{G2}, it is natural to
 see the set $\go{bb}$
of normalized  block sequences of a given basis as a
closed subset of ${\go X}^w$ equipped with the product
of the norm topology; as expected,
$\approx$ is Borel in $\go{bb}^2$.
The associated canonical injection from $2^w$ (resp. $\go{bb}$) into $\go A$
is Borel.

Going back to the notion of isomorphism, we see that it is indeed analytic:
For $\go X, \go Y \in \go A$ we have $\go X\iso \go Y$ iff
$$\e \pmb x \; \e \pmb y \; (\pmb x \approx \pmb y \; \land\; \a n \; x_n \in
\go X\; \land\; \a n \; y_n \in \go Y \;\land\; \a n (\ku U_n \cap \go X\neq
\tom \to \e m\; x_m \in \ku U_n) $$
$$ \land\; \a n (\ku U_n \cap \go Y\neq \tom \to \e m \; y_m \in \ku U_n))$$

Here $\mgdv \ku U_n\mgdh_{n \in \N}$ is a basis for the topology on $C(\ca)$.

Our first result concerns
equivalence of subsequences (resp. block-bases) of a given basis.

 \section {Counting the number of nonequivalent\\ blockbases.}

 Until now we have only been looking at Cantor`s concept of
 cardinality, but there is also a newer and finer one steming from
 descriptive set theory. It allows us to distinguish between different
 levels of $\ca$ according to their complexity. For the definition:

 \begin {defi} Let $E \subset \go X^2$ and $F \subset \go Y^2$ be
equivalence relations
 on standard Borel spaces $\go X$ and $\go Y$. We say that $E$ is
 Borel reducible to $F$ ($E \le_{B} F$) if
there is a Borel function $\phi
: \go X
 \longrightarrow \go Y$ such that $\a x,y \in \go X \quad xEy
 \leftrightarrow \phi (x)F \phi (y)$.
$E$ is Borel bireducible with $F (E \sim_{B} F)$ if both $E \le_{B} F$
and $F \le_{B} E$.

 \end {defi}

 The definition has not much interest unless $\go X$ and $\go Y$ are
 uncountable and $E$ and $F$ have uncountably many classes.
 Furthermore it is usually supposed that the equivalence relations
 are of some bounded complexity, eg. Borel or analytic.

 Being given some normalized basic sequence  $\{e_i\}_{i \in \N}$ in
 a separable Banach space $\go X$, we can look at the set of its subsequences
 as
 a subspace of $\ca$, just take away the countable set $FIN$ and
 identify a subsequence with its set of indices. $\ca \backslash FIN$
 is still a Polish space (thoughnow noncompact) under the usual topology.

  Note that equivalence of such normalized basic sequences, $\approx$
induces a Borel
equivalence relation on
$\ca
 \backslash FIN$.

 An important measure of complexity is the following equivalence
 relation $E_{0}$, it is the minimum (up to $\sim_{B}$) Borel equivalence
relation
 $\leq_{B}$ above identity on Polish spaces and is defined on $\ca$ as follows:
 $$\alpha E_{0} \beta :\equiv \e n \in \N \quad \a m \geq n \quad
 \alpha_{m}=\beta_{m}$$

 \begin {prop}Let  $\go X$ be a Banach space with
a basis $\{e_i\}_{i \in \N}$.

 \begin {itemize}
 \item Either  $\{e_i\}_{i \in \N}$ is perfectly homogeneous, i.e.
 equivalent to all of its normalized block basic sequences (and
 therefore equivalent to the standard unit basis in some $l_{p},
 1\leq p< \infty$ or $c_{0}$), or $E_{0}\leq_{B}
 \approx\upharpoonright_\go{bb}$;
\item Either $\{e_i\}_{i \in \N}$ is subsymmetric, i.e. equivalent to
 all of its subsequences,  or $E_{0}\leq_B \sim$.

\end{itemize}
\end{prop}

\pf
We show only the first part as the proof of the second is
essentially the same.

Assume that $\go X$ has at least two non-equivalent normalized block sequences,
$(x_i)$ and $(y_i)$. Then $(x_i)$ and $(y_i)$ are not $2$-equivalent,
 so there
exists $I_1$ an interval of integers such that $(x_i)_{i \in I_1}$ and
$(y_i)_{i \in I_1}$ are not $2$-equivalent.
Also for any $k$, $(x_i)_{i>k}$ and $(y_i)_{i>k}$ are not equivalent so
there exists
$k <I_2 \subset \N$ such that $(x_i)_{i \in I_2}$ and
$(y_i)_{i \in I_2}$ are not $4$-equivalent.
 Without loss of generality,
we may assume that $supp(x_i, i \in I_2) \cup supp(y_i, i \in I_2)
> supp(x_i, i \in I_1) \cup supp(y_i, i \in I_1)$
Repeating
this, we get a sequence of successive intervals $I_k$ such that $(x_i)_{i
\in I_k}$ and
$(y_i)_{i \in I_k}$ are not $2^k$-equivalent.

For any $\alpha \in \ca$, let $f(\alpha)$  be the unique
block
basis defined by the set of vectors $\{f_i | 
i \in \cup_{j \in \N}I_j\}$, with
$f_i=x_i$ if $i \in  I_j$ with $\alpha(j)=0$, $f_i=y_i$ 
if $i \in I_j$ with $\alpha(j)=1$; clearly,
$f(\alpha)$ is equivalent to $f(\beta)$ iff $\alpha E_0 \beta$
 and the map is
Borel (even continuous).\pff

Let us mention that in the case where $E_0$ Borel reduces to an analytic
equivalence relation this latter can have no analytic section. That is, there
can be no analytic set intersecting every equivalence class in
exactly one point. So in particular the above result says that if there is
more that one class, then the relation is rather complicated, and cannot in
fact be classified in a Borel manner by real numbers.
Also, we refer to the forthcoming \cite{FR} for discussion
about isomorphism on $\go{bb}$.

\

\subsection{Subsequences of a basis}

\subsubsection{Classical Ramsey type results.}

We first mention that classical Ramsey type results can be proved in the
case where
there are countably many classes of isomorphism.

\begin{lemme} For any $n$, there exist a constant $c(n)=1+n(2^{n+1}+1)^2$,
such that for any Banach space $\go X$, all $n$-codimensional subspaces of
$\go X$
are
$c(n)$ isomorphic.
\end{lemme}
\pf
We are not interested in finding the best constant. We prove that for any Banach
space $\go X$, any $n$-codimensional subspace
$\go Z$ of $\go X$, $d_{BM}(\go X,{\go Z} \oplus_1
l_2^n) \leq \sqrt{n}(2^{n+1}+1)$. The result follows: for any two
$n$-codimensional subspaces
$\go H$ and
$\go H'$ of a Banach space, there exists
 $\go W$ $n$-codimensional in $\go H$ and
$\go H'$, so $\go H$ (resp.
$\go H'$) is
$\sqrt{n}(2^{n+1}+1)+\eps$-isomorphic to
$\go W
\oplus_1 l_2^n$ for all $\eps>0$, and $d_{BM}(\go H,\go H') \leq
n(2^{n+1}+1)^2$.

So let $\go X$ be a Banach space and $\go Z$ be $n$-codimensional in $\go X$. By
induction, there
exists a projection $P$ on $\go Z$ of norm smaller than $2^n+\eps$. Let
$F=(Id-P)(\go X)$. Then
using $P$ one shows that
$\go X=\go Z \oplus F$
is $2^{n+1}+1+2\eps$-isomorphic to $\go Z \oplus_1 F$. By classical results, $F$
is
$\sqrt{n}$-isomorphic to $l_2^n$,
and it follows that $\go X$ is $\sqrt{n}(2^{n+1}+1+2\eps)$-isomorphic to $\go Z
\oplus_1 l_2^n$.
\pff

To spark
confusion we will identify several different objects; namely the
space $[e_{i}]_{i\in A}$ for some subset $A\subset \N$ with the
characteristic function $\chi_{A}$ seen as a point in the Cantor
space $2^\omega$ which we again simply identify with the subset A of
$\N$. The subspaces spanned by subsequences are therefore equipped with
the compact Polish topology inherited from $2^\omega$.
We denote by $A\cong B$ the fact that the corresponding Banach spaces are
linearly isomorphic. If we see this relation as an equivalence
relation between the points in the Effros Borel space of closed linear
subspaces of $C(\ca)$ it is analytic, and furthermore the function associating
to $A\subset\N$ the space (or point) in the Effros Borel space is
Borel. So, as noted at the end of first section, we have an induced
analytic equivalence
relation,
denoted by $\cong$, on $2^\omega$. For $K \geq 1$, write $A \cong_K B$ to
mean that
the corresponding spaces are $K$-isomorphic, so that
$\cong$ is the union of the $\cong_K, K \in \N$.
Recall that for $u=<u(0),\ldots,u(k)> \in 2^{<\omega}$, $N_u$ denotes the
basic open set
of all $A \in \ca$ such that $i \in A$ iff $u(i)=1$, $\forall i \leq k$.

\begin{prop}
Let $\go X$ be a Banach space with a basis $\{e_i\}_{i \in \N}$. If $\go X$ has countably many
classes
of isomorphism
generated by subsequences of the basis  then there exists $K$ and
$A \subset \N$ such that $[e_i]_{i \in B} \sim_K [e_i]_{i \in A}$ for all
$B \subset A$.
\end{prop}

\pf
In this proof, we think of sets $A$ as increasing sequences of integers.
For each $n$, let $E_n \subset \N$ be a
representative in the $n$-th $\cong$ class.
Say that $A$ belongs to ${\cal A}_{m,n}$ if $A \cong_m E_n$; fix furthermore
a bijection $k\mapsto (m_k,n_k)$ between $\N$ and $\N^2$.

Assume a Banach space $\go X$ contradicts the proposition
 and apply the infinite Ramsey theorem.
No subset
of $\N$ has all its subsets in ${\cal A}_{m_0,n_0}$, so there
exists $A_0 \subset \N$ so that no subset of $A_0$ is in ${\cal A}_ {m_0,n_0}$.
Let $k_0$ be
 an integer in $A_0$.
No subset $A$ of $A_0$ after $k_0$ is such that $\{k_0\} \cup A$ has all
its subsets
 in ${\cal A}_{m_1,n_1}$: otherwise, any subset $\{b_0<b_1<\ldots\}$ of $A$,
spans a space which is $c(1)$-isomorphic
to the span of $\{k_0<b_1<\ldots\}$, so is $c(1)m_0$ isomorphic to $E_{n_0}$, a
contradiction. So some subset $A_1$ of $A_0$ after $k_0$ is such that $\{k_0\}
\cup A$
 belongs to ${\cal A}_{m_1,n_1}$ for no $A \subset A_1$. Choose $k_1$ in
$A_1$. Repeating this, define a decreasing sequence $A_n$ and an increasing
sequence $k_n$
with $k_i$ in $A_i$, to get a  set $K=\{k_p, p \in \N\}$ belonging to no
${\cal A}_{m,n}$.
This should be sufficient to convince you. \pff

The example of  $l_1 \oplus l_2$ proves that the isomorphism class in this
proposition
need not contain $\go X$, nor contain "most" subspaces of $\go X$, in fact, we
shall see that
the class of $l_1$ (resp. $l_2$) is meager in the standard topology
associated to the space.

\section {Banach spaces with unconditional basis\\
isomorphic to their squares.}

Let $\go X$ be a Banach space with a basis $\{e_i\}_{i \in \N}$. We
will study the linear-isomorphism classes of spaces spanned by subsequences
(finite and
infinite) of this basis by using
 Baire category. \\

As the finite subsets of $\N$ are fully characterized up to $\iso$ by their
 cardinality we will often forget about them, therefore we want a name for
 them, so we remember what we are forgetting. Consequently let $FIN$ be the
Frechet
 ideal of finite subsets of $\N$.

Also other relations on $2^\omega$ will be useful. First note that the
Cantor space is an Abelian Polish group under the action of symmetric
difference
$\triangle$, as such the identity element is $\emptyset$ and each
element is its own inverse. Therefore not only $\cong$ but also the
relations :
$$A\sim_{1}B : \equiv A \cong \complement B$$
$$ A\sim_{2}B : \equiv \complement A \cong \complement B$$
are analytic in $2^\omega$.

\begin{thm}{(Kuratowski/Mycielski)}
 Let $\mathfrak {X}$ be a perfect Polish space, and
R be a relation on $\go X$ meager in $\go X^{2}$. Then there exists a
homeomorphic copy $\ku C$ of the Cantor space such that $\forall x,y \in
\ku C$  with $x\ne y$ we have $\neg xRy.$
\end{thm}

Suppose now that $\cong$, $\sim_{1}$ and $\sim_{2}$ all are meager in
the product $2^\omega \times 2^\omega$, then so is their union and the
theorem above gives us a Cantor set $\ku C$ avoiding the three relations.
Now by taking away the countable set (possibly finite or empty) of finite
sets in $\ku C$ and going to a further subset, we can suppose that
$\ku C\cap FIN =\emptyset$.

Following Kalton we say that the basis $\{e_{i}\}_{i \in \N}$
is
countably primary if there is a countable list $E_{0}, E_{1}, E_{2},
\ldots$ of Banach spaces such that if $A \subset \N$ then for some $n$ either
$[e_{i}]_{i \in A} \cong E_{n}$ or $[e_{i}]_{i \in \complement A}
\cong E_{n}$.
It is now easy to see that our space cannot be countably primary, for
take two subsets of $\N$, $A \ne B, A,B \in \ku C$, such that they have the
same
$E_{n}$ associated, i.e. either $A \cong E_{n} \cong B$, $A\cong E_{n}
\cong \complement B$, $\complement A \cong E_{n} \cong B$ or $\complement A
\cong E_{n} \cong \complement B$.
It should be evident for the reader taking the pains of unraveling
our definitions that we get a minor contradiction here.

Going in the other direction we wish to see what we get from the fact
that some of the three relations are nonmeager in the product. First
we note that a Fubini type theorem is true also in Baire category:

\begin{thm}{(Kuratowski/Ulam)}
   Let $\go Y$ be a Polish space and $\ku D$ a subset of the square having
    the Baire property, then \\
    $\ku D$ is nonmeager $\equiv \exists^{*}x \exists^{*}y (x,y) \in
    \ku D
    \equiv \exists^{*}y \exists^{*}x (x,y) \in \ku D$
    \end{thm}

    Here $\exists^{*}x P(x)$ signifies the existence of a nonmeager
    set of $x$ such that $P(x)$.
    Applying this to the above relations we get the following:

    $$\cong \textrm{nonmeager in}\ 2^\omega \times 2^\omega \Rightarrow
    \exists A \in 2^\omega \exists^{*}B \in 2^\omega A \cong B$$

    $$\sim_{1} \textrm{nonmeager in}\ 2^\omega \times 2^\omega \Rightarrow
\exists B \in 2^\omega \exists^{*}A \in 2^\omega A \cong \complement
B $$ $$\Rightarrow
\exists C \in 2^\omega \exists^{*}A \in 2^\omega A \cong C$$

$$\sim_{2} \textrm{nonmeager in} \ 2^\omega \times 2^\omega \Rightarrow
\exists A \in 2^\omega \exists^{*}B \in 2^\omega \complement A \cong \complement
B $$ $$\Rightarrow \exists C \in 2^\omega \exists^{*}D \in 2^\omega C \cong
D  $$

So all of the three cases gives us a nonmeager isomorphism class
$\mathcal {A} $
in $2^\omega$. Fix $B$ in $\mathcal {A}$ and for all $M \in \N$ denote by
 $\mathcal {A}_M $ the set $\{A: A \cong_M B\}$. Then for
some $M$ large enough in $\N$, $\mathcal {A}_M $ is nonmeager as well, and we
intend to show that $\mathcal {A}_{K}$ for some $K \geq M$ (and therefore
$\mathcal{A}$) is actually
residual in the Cantor
space.

Note that being a section of an analytic set ($\mathcal
{A}= (\cong_M)_{B}$) $\mathcal{A}_M$ is itself
analytic and has therefore the property of Baire. So as it is nonmeager
it must be residual in an open set $\mathcal {U} \subset 2^\omega$ and
by going to a smaller open set we can suppose that $\mathcal {U}$ is
on the form $N_{s}$ for some finite sequence $s \in 2^{<\omega}$.
 There are for us two
interesting features of $s$, its length
and its
cardinality. By the length, $|s|$, we denote its length as a
sequence and by the cardinality, $\overline s$, is denoted the number
of $1$'s appearing in the sequence.
We prove that $\mathcal {A}_K$ is residual for $K=Mc(2|s|)$.

Otherwise, let $t \in 2^{<\omega}$ be such that $\mathcal {A}_K$
is meager in $N_t$. Without loss of generality, assume $|t| \geq |s|$
and write $t=u^\frown v$, with $|u|=|s|$. Then let
$t'=t^\frown s=u^\frown v^\frown s$ and
$s'=s^\frown v^\frown u$.
So $|t'|=|s'|$, $\overline t'=\overline s'$,
and $\mathcal {A}_M$ is residual in $N_{s'}$
while $\mathcal {A}_K$ is meager in $N_{t'}$.

We have now a natural
homeomorphism $\phi$ between the clopen sets $N_{s'}$ and $N_{t'}$; simply for
an $A \in N_{s'}$ change the beginning from $s'$ to $t'$, i.e. $\phi
({s'}^\frown \alpha)= {t'}^\frown \alpha$, where $\chi_{A}={s'}^\frown \alpha$.
 By construction,
${s'}^\frown \alpha$ and ${t'}^\frown \alpha$ code subspaces of same
codimension
$2|s|-\overline{s}-\overline{u} \leq 2|s|$
of the space coded by $(1^{|s|})^\frown v^\frown (1^{|s|})^\frown \alpha$
(here $1^n$ denotes the length $n$ sequence of $1$'s).
By Lemma 3, it follows that $\phi(A) \cong_{c(2|s)} A$.
As $\mathcal {A}_M$ is residual in $N_{s'}$, it follows that
$\mathcal {A}_{c(2|s|)M} \supset \phi(\mathcal {A}_M)$ is residual in $N_{t'}$,
 a contradiction.

We now need a standard compactness result from descriptive set theory:

\begin {lemme} If $\ku {G}$ is a residual subset of $\ca$, then there
exists a partition \\ $A_{0}, A_{1}, A_{2}, \ldots$ of $\N$ and subsets
$B_{i} \subset A_{i}, i \in  \N$ such that for any set $E \subset \N$
if $\e i \in \N \quad E \cap A_{i}= B_{i}$ then $E \in \ku G$
\end {lemme}
\pf For $\ku D$ an open dense set in $\ca$ and $n\in \N$
there is $s \in 2^{<\omega}$ such that for any $t \in 2^n$ we have
$N_{t^\frown s} \subset \ku D$.

This because we can ennumerate $2^n$ as $\{
t_{1},t_{2},\ldots,t_{2^n}\}$ and so due to the density and the fact
that $\ku D$ is open there is some $s_{1} \in 2^{<\omega}$ with
$N_{t_{1}^\frown
s^ {}_{1}} \subset \ku D$.

Find $s_{2} \in 2^{<\omega}$ with $N_{t_{2} ^\frown s_{1}
^\frown s_{2}} \subset \ku D$.

Again find $s_{3} \in 2^{<\omega}$ with $N_{t_{3} ^\frown s_{1}
^\frown s_{2} ^\frown s_{3}} \subset \ku D$, etc.

Now put $s=s_{1}^\frown s_{2}^\frown \ldots ^\frown s_{2^n}$ and we
 have our result.

 Suppose that $\ku G \supset \snit_{i \in \N} \ku D_{i}$ for open, dense
 sets $\ku D_{i}$.
 \begin {itemize}
 \item Take $s(0,0)$ such that $N_{s(0,0)} \subset \ku D_{0}$,\\ put
 $n(1,0)=|s(0,0)|$.
 \item Take $s(1,0)$ such that $N_{t ^\frown s(1,0)} \subset \ku D_{0}, \a
 t \in 2^{n(1,0)}$,\\ put
 $n(0,1)=n(1,0)+|s(1,0)|$.
\item Take $s(0,1)$ such that $N_{t ^\frown s(0,1)} \subset \ku D_{1}, \a
 t \in 2^{n(0,1)}$, \\put
 $n(2,0)=n(0,1)+|s(0,1)|$.
\item Take $s(2,0)$ such that $N_{t ^\frown s(2,0)} \subset \ku D_{0}, \a
 t \in 2^{n(2,0)}$,\\ put
 $n(1,1)=n(2,0)+|s(2,0)|$.
\item Take $s(1,1)$ such that $N_{t ^\frown s(1,1)} \subset \ku D_{1}, \a
 t \in 2^{n(1,1)}$, \\put
 $n(0,2)=n(1,1)+|s(1,1)|$.
 \item Take $s(0,2)$ such that $N_{t ^\frown s(0,2)} \subset \ku D_{2}, \a
 t \in 2^{n(0,2)}$, \\put
 $n(3,0)=n(0,2)+|s(0,2)|$,
 \end{itemize}
 etc.
 \begin{itemize}
 \item $A_{i}:=\for_{k \in \N} \intv n(i,k),n(i,k)+|s(i,k)|\intv$
 \item $B_{i}:=\for_{k \in \N} \mgdv n(i,k)+m\;\big\vert\;  (s(i,k))(m)=1
\mgdh$\pff
\end {itemize}

Letting $A'_{0}:=A_{0}, A'_{1}:= \for _{i=1}^{\infty} A_{i},
B'_{0}:=B_{0}, B'_{1}:= \for _{i=1}^{\infty} B_{i}$ we get:
If $\mathcal {G}$ is a residual set in $2^\omega$, then there is a
partition $\N=A'_{0}\cup A'_{1}, A'_{0} \cap A'_{1}= \emptyset$ and sets
$B'_{i} \subset A'_{i}, i \in {0,1}$, such that for any $D \subset \N$,
if either $D\cap A'_{0}=B'_{0}$ or $D\cap A'_{1}=B'_{1}$ we have $D \in
\mathcal {G}$.

We shall now assume the basis is unconditional.
Apply this to $\mathcal {A}$  and we see in particular that $B_{0}$,
$B_{1}$ and $B_{0} \cup B_{1} \in \mathcal {A}$.
Now as $\{e_i\}_{i \in \N}$ is unconditional, we have for any pair of
disjoint sets $C, D \subset \N$ that $[e_i]_{i \in C} \oplus
[e_i]_{i \in D} \cong [e_i]_{i \in C \cup D}$. Moreover
$\mathcal {A}$ being residual and $\complement$ a homeomorphism of
$2^\omega$ there is some $C\subset \N$ with $C, \complement C \in
\mathcal {A}$. So again abusing notation we calculate:
$$\N \cong C \cup \complement C \cong C \oplus \complement C \cong
B_{0} \oplus B_{0} \cong B_{0} \oplus B_{1} \cong B_{0} \cup B_{1}
\cong B_{0}.$$
 So $\N$ belongs to the residual class $\mathcal {A}$
 and
$[e_i]_{i \in \N}$ is isomorphic to
$[e_i]_{i
\in B_{0}}$ which
 is isomorphic to its square.

Without loss of generality assume $\N$ belongs to
${\cal A}_K$. Let $c$ be the
constant of unconditionality of the basis. We denote by $\oplus_1$
the $l_1$-sum of Banach spaces.

 Again take some arbitrary subset $D \subset \N$, we have due to the
 fact that $A_{0},A_{1}$ partition $\N$:
 $$\N \oplus_1 D
 \cong_K (B_{0} \cup  B_{1}) \oplus_1 ((D \cap A_{1}) \cup
 (D\cap A_{0}))
 \cong_{2C} B_{0} \oplus_1 B_{1} \oplus_1 (D \cap A_{1}) \oplus_1
 (D\cap A_{0})
 $$
 $$\cong_{2C} (B_{0} \cup (D \cap A_{1})) \oplus_1 (B_{1}
 \cup (D \cap A_{0})) \cong_K B_0 \oplus_1 B_1 \cong_K \N$$

So spaces of the form $[e_i]_{i \in \N} \oplus_1 [e_i]_{i \in D}$ for any
 $D\subset \N$ are $4C^2 K^3$-isomorphic to $[e_i]_{i \in \N}$; and
 in particular $[e_i]_{i \in \N}$ is isomorphic to its hyperplanes.

 Using also the complete version of the above lemma, we find \\
 $\for _{i \in \N} B_{i}, B_{0}, B_{1}, B_{2}, \ldots \in \ku {A}$.
 But $[e_{i}]_{i \in B_{0}}, [e_{i}]_{i \in B_{1}}, [e_{i}]_{i \in
 B_{2}}, \ldots$ gives (due to the unconditionality of the basis) an
 unconditional Schauder decomposition of $[e_{i}]_{i \in \for _{i \in
 \N} B_{i}}$ and this latter is again isomorphic to $[e_{i}]_{i \in \N}$.

 Summing up we have arrived at the following:

 \begin {thm}{} Let $\go X$ be a Banach space with an unconditional basis
$\{e_i\}_{i \in
\N}$, then either:

(1) There exists a perfect set $\ku P \subset 2^\omega$ of infinite
subsets of $\N$ such that for any two distinct $A,B \in \ku P$ we have that
$$[e_i]_{i \in A} \ncong [e_i]_{i \in B}$$
 $$[e_i]_{i \in A} \ncong [e_i]_{i \in \complement B}$$
 $$[e_i]_{i \in \complement A} \ncong [e_i]_{i \in \complement B}$$

 or

 (2) For $A$ in a residual subset of $2^\omega$, $[e_i]_{i \in A}$
is isomorphic to $\go X$; $\go X$ is isomorphic to
its hyperplanes, to its
 square, uniformly isomorphic to $\go X \oplus [e_i]_{i \in D}$
for any
 $D\subset \N$ and to a denumerable unconditional Schauder decomposition into
 uniformly isomorphic copies of itself, $\go X \iso (\sum _{k \in
 \N} \oplus ([e_i]_{i \in B_{k}}))$.
 \end {thm}

 As we noticed earlier, (1) in our dichotomy implies that the least
 cardinal $\kappa$ such that the basis is $\kappa$-primary (with the
obvious definition) is $\ca$,
 that is the trivial one.

The theorem improves on an earlier result by Kalton which gives (2) in case
the unconditional basis is countably primary. His proof like the above used
very little from Banach space theory, instead his setting was measure and
probability theory \cite{kal}.

As is easily seen from the proof, (1) can be strenghtened considerably
and to our purposes, in fact one can get a perfect set avoiding any
countable list of relations of the form
$$[e_i]_{i \in A} \cong \phi (B)$$
where $\phi$ is a Borel function from $\ca$ to the Effros Borel space
of closed linear subspaces of $C(\ca)$, that is the canonical space
of separable Banach spaces. For example one
could use Borel functions $\psi$ from $\ca$ to the space $\go{bb}$ of normalized
block sequences of some given basis, and one would avoid
$$[e_i]_{i \in A} \cong [\psi (B)].$$

Eg. we can force
$$[e_i]_{i \in A} \ncong ([e_i]_{i \in B})^2$$
$$[e_i]_{i \in A} \ncong ([e_i]_{i \in B})^3$$
$$[e_i]_{i \in A} \ncong ([e_i]_{i \in B})^4$$
and so on.
Or
$$[e_i]_{i \in A} \ncong c_{0}([e_i]_{i \in B})$$
or whatever construction from B being reasonably explicit.

\begin {ex}
There is a certain sense in which the above result is
optimal. For we might like to try to show that not only
is some isomorphism class residual in $\ca$ but that it is
all of $\ca \backslash FIN$. But this is easily seen to be false,
for take the following basis for $l_1\oplus l_2$:
$$\|a_0e_0+a_1e_1+\ldots+a_{2n+1}e_{2n+1}\|:=|a_0|+|a_2|+\ldots+|a_{2n}|+
\sqrt{a_1^2+a_3^2+\ldots+a_{2n+1}^2}$$
then there are exactly three isomorphism classes: $l_1$, $l_2$, and $l_1\oplus
l_2$ the first two being meager in $\ca$:
Cause if $A \subset \N$, $A$ infinite, contains infinitely many even and
odd numbers then
$\intv e_i \inth_{i \in A} \iso l_1 \oplus l_2$, if it only contains
finitely many even numbers then $\intv e_i \inth_{i \in A} \iso l_2$, and
only finitely many odd numbers then $\intv e_i \inth_{i \in A} \iso l_1$.
\end {ex}

\begin {ex}
If we take the standard Haar basis for some $L_p ([0,1])$, $1<p<\infty$
then this is unconditional and in fact the only two spaces spanned by
subsequences are $l_p$ and  $L_p ([0,1])$; so there are bases inducing
exactly two isomorphism classes in $\ca \backslash FIN$. \cite{lt2} Thm.
2.d.10.
\end {ex}

\begin{ex}Tsirelson's space.
\end{ex}

 We take a look at the standard unit vector basis for the Tsirelson space.
It has the following properties (here $\N^*=\N \backslash \{0\}$):

\begin {itemize}
\item Two subsequences $\{t_{k_{i}}\}_{i \in \N^{*}}$ and $\{ t_{l_{i}} \}_{i
\in \N^{*} }$
are equivalent iff $$[t_{k_{i}}]_{i \in \N^{*}} \iso [t_{l_{i}}]_{i \in
\N^{*}}$$
\item Two subsequences $\{t_{k_{i}}\}_{i \in \N^{*}}$ and $\{t_{l_{i}}\}_{i
\in \N^{*}}$
are equivalent iff
$$sup_{i \in \N^{*}} \big ( \|I_{i}\|,\|J_{i}\|  \big ) < \infty$$
where
$$I_{i}:\intv t_{l_{(l_{n})}} \; \big| \; k_{i-1}<l_{n}\leq
k_{i} \inth \hookrightarrow l_{1}$$
and
$$J_{i}:\intv t_{k_{(k_{n})}} \; \big| \; l_{i-1}<k_{n}\leq
l_{i} \inth \hookrightarrow l_{1}$$
are the formal identities ($k_{0}:=l_{0}:=0$).

\item$ \{t_i\}_{i\in \N^*}$ and $\{t_{k_{i}}\}_{i \in \N^{*}}$ are
equivalent iff the function $k:\N^* \longrightarrow \N^*$ is majorized by a
primitive recursive function.
\end{itemize}

 We can see subsets $A \subset \N^*$ as strictly increasing functions
$a:\N^* \longrightarrow \N^*$, simply let $a$ enummerate $A$ in the usual
order. In the same way, strictly increasing functions can be seen as
infinite subsets of $\N^*$.
Now the relation of $b$ majorizing $a$  is in fact closed in $(\ca
\backslash FIN)^2$ and has closed, nowhere dense sections:

$$a \leq b \equiv \a n\; a(n)\leq b(n) \equiv $$
 $$\a n \; \a m,k \Intv \intv \a m_1<m_2<\ldots<m_n<m  \;   \e i B(m_i)=0
\og $$
$$\e m_1<m_2<\ldots<m_n=m
 \a i B(m_i)=1 \og$$
$$ \a k_1<k_2<\ldots<k_n <k \;  \e i A(k_i)=0  \og $$
$$\e k_1<k_2<\ldots<k_n=k   \;
 \a i A(k_i)=1 \inth \to k\leq m \Inth$$

So it's closed and it is easily seen that no function $a$ can be such that
it majorizes all functions belonging to some basic open set $N_s$ (on the
other hand $n\mapsto n$ minorizes all functions), whence the sections (in
one of the coordinates) have empty interior.

The set $\mgdv A \subset \N^* \del A\iso \N^* \mgdh$ is the countable union
of closed nowhere dense sets, hence meager.

Now take any $A \not \iso \N^*, A=\mgdv a_n\mgdh_{n\in \N^*}$. As it
generates a sequence non equivalent with the full basis there are disjoint
nonempty intervals
$$\intv r_1,s_1\intv, \intv r_2,s_2\intv,\ldots \subset \complement A$$
such that for the formal identities
$$I_i:\Intv t_n \Del n\in \intv r_i,s_i\inth \Inth \hookrightarrow l_1$$
we have $\|I_i\| > i$.
So if $B\iso A$ then $\e i \; \a j \geq i \; \intv r_j,s_j \intv \not
\subset B$, which again easily is seen to be the countable disjunction of
closed, nowhere dense conditions.

Therefore every isomorphism class for the Tsirelson space is meager in $\ca
\backslash FIN$. However, using the functions $n\mapsto 2n$ and $n\mapsto
2n+1$, it is seen that the space is isomorphic to its square.

\medskip

Though the notion of `complement' of a set is not absolute, ie. it depends
on the ambient space, $\iso$ is so.
This is to say that $A\iso B$ does not depend on whether $A$ and $B$ are
seen as subsets of $\N$ or of  any other $C \subset \N, A,B\subset C$.

Denote by $2^A$ the closed set $\mgdv \chi_B \in \ca \del B\subset A\mgdh$.
If $A$ is infinite then this set is homeomorphic to $\ca$ and we can use
our preceeding arguments thereon.

So if for some $A$ the restriction of $\iso$ to $2^A$ is meager in $2^A$
there are $\ca$ classes of isomorphism generated by subsequences of $A$,
hence also of $\N$.

If  not there is some residual in $2^A$ isomorphism class and we get as
before that
 $[e_i]_{i \in A}$ is isomorphic to its hyperplanes, to its
 square, to $[e_i]_{i \in A} \oplus [e_i]_{i \in D}$ for any
 $D\subset A$ and to a denummerable Schauder decomposition into
 isomorphic copies of itself, $[e_i]_{i \in A} \iso (\sum _{k \in
 \N} \oplus ([e_i]_{i \in B_{k}}))$.

There is a priori no control on the uniformity of the isomorphisms between
each space $[e_i]_{i \in A}$ and its square. Notice also that
the existence of an isomorphism $A \oplus D \iso A$ for $D$ infinite subset
of $A$
 is straightforward
from the assumption that $A \iso A \oplus A$ for all infinite $A \subset \N$.
For then $A \oplus D \iso (A \setminus D) \oplus D \oplus D
\iso (A \setminus D) \oplus D \iso A$, for all infinite $D \subset A$.

The diligent reader is invited to amuse himself in applying the above proof
to the cases of Lipschitz homeomorphism, uniform homeomorphism and
permutative equivalence of bases (all three are analytic equivalence
relations).

\begin {thm} Let $\go X$ be a Banach space with an unconditional basis
$\{e_i\}_{i \in \N}$. Then either:

(1) There exists a perfect set $\ku P \subset 2^\omega$ of infinite
subsets of $\N$ such that for any two distinct $A,B \in \ku P$ we have that
$$[e_i]_{i \in A} \ncong [e_i]_{i \in B}$$

 or

 (2)For any infinite subset $A\subset \N$:
 $[e_i]_{i \in A}$ is isomorphic to its hyperplanes, to its
 square,  and to a denumerable Schauder
decomposition into uniformly
 isomorphic copies of itself, $[e_i]_{i \in A} \iso (\sum _{k \in
 \N} \oplus ([e_i]_{i \in B_{k}}))$.
 \end {thm}

Recall that a Banach space is \em{complementably minimal} if it embeds
complementably in any of its subspaces. Using Pelczynski's decomposition
method and the above result, one proves:

\begin {cor} Let $\go X$ be a complementably minimal Banach space with
an unconditional basis $\{e_i, i\in\N\}$. Then

(1)  There exists a perfect set $\ku P \subset 2^\omega$ of infinite
subsets of $\N$ such that for any two distinct $A,B \in \ku P$ we have that
$$[e_i]_{i \in A} \ncong [e_i]_{i \in B}$$

or

(2) For any infinite subset $A \subset \N$,
$[e_i]_{i \in A}$ is isomorphic to $\go X$.

\end {cor}

\section {The number of nonisomorphic subspaces of \\hereditarily indecomposable
 Banach spaces.}

This result was proved by the second-named author before the authors's
collaboration began.

 Let now $\go X$ be some separable hereditarily indecomposable (H.I.) Banach
 space. This means that no (closed, infinite dimensional) subspace
 of $\go X$ can be written as a direct sum of two closed infinite
 dimensional subspaces. In this section `space' will always refer to
 closed infinite dimensional subspaces of $\go X$.
It follows clearly from the H.I. property that $\go X$
contains no unconditional basic sequence; and in fact, by Gowers's dichotomy
theorem, every Banach space contains either a H.I. subspace or
a subspace with an unconditional basis.
Moreover Gowers and Maurey proved that a H.I. space is isomorphic to no proper
subspace
(and as the H.I. property is hereditary, this is also true of any
subspace of $\go X$).

 Now since we have the first property we cannot hope to use the
 above theorem to conclude something about the number of
 nonisomorphic subspaces of $\go X$, but it is  still possible to use the fact
 that $\go X$ must contain some basic sequence $\{e_i\}_{i \in \N}$.
Again we look at the subspaces spanned by subsequences as points in
 $\ca$.

 According to Dedekind every real is a set of
 rational numbers $r= \mgdv q \in \Q | q<r \mgdh$, but following Cantor
 the set of rational numbers is the
 same as the set of natural numbers. So every real $r$ is a set
 $A_{r}$ of naturel numbers such that $r<s \equiv A_{r} \subsetneq A_{s}$.
 Now again confusing Banach spaces with subsets of $\N$, reals become
 for us a subspace $\go B_{r}$ and the relation $r<s$ is simply strict
 inclusion $\go B_{r} \subsetneq \go B_{s}$.
 All of the identifying functions are of course (and evidently) Borel.
 So we have $\ca$ nonisomorphic subspaces of $\go X$.

 This could also
 have been seen using the arguments from the preceeding section, for
 suppose that some isomorphism class $\ku A$ was residual in $\ca$ then
 as $\{ 0 \} \triangle (\cdot)$ is a homeomorphism of $\ca$
 $$\ku A \cap \mgdv \{ 0 \} \triangle A \; | \;  A \in \ku A \mgdh
 \neq \tom$$
 So there is some $B \subset \N, 0 \not \in B$ with $B, \{ 0 \} \
 B \in \ku A$, ie. some subspace of $\go X$ isomorphic to a hyperplane.
 But this cannot be the case in an H.I. space, so $\iso$ must be meager
 in $\ca$ and the Kuratowski/Mycielski result takes care of the rest.

 \begin {prop}
 Any H.I. Banach space contains $\ca$ pairwise non isomorphic subspaces.
 \end {prop}

\begin{cor}
Any Banach space contains $2^{w}$ pairwise non isomorphic subspaces or
is saturated with subspaces isomorphic to their
squares.
\end{cor}

\

\

Equipe d'Analyse,

Couloir 46-0, Boite 186,

Universit\'{e} Paris 6,

4, place Jussieu,

75252 Paris Cedex 05,

FRANCE.

E-mail: ferenczi@ccr.jussieu.fr,

rosendal@ccr.jussieu.fr

\begin{thebibliography}{AAA}

\bibitem{FR} V. FERENCZI, C. ROSENDAL, {\em Countably homogeneous spaces},
in preparation.


\bibitem{G1} W.T. GOWERS, {\em A new dichotomy for Banach spaces},
 Geometric and Functional Analysis {\bf 6} (1996), 6, 1083-1093.

\bibitem{G2} W.T. GOWERS, {\em An infinite Ramsey Theorem and some
Banach-space
dichotomies}, preprint.

\bibitem{GM} W.T. GOWERS and B. MAUREY, {\em The unconditional basic
sequence problem}, J.Amer.Math.Soc. {\bf 6} (1993), 851-874.

\bibitem{kal}N. KALTON: A remark on Banach spaces isomorphic to their squares,
\emph{Contemporary Mathematics, Volume $\fed{232}$, (1999)}

\bibitem{ke}A. KECHRIS: Classical Descriptive set theory, \emph{Springer,
New York (1995)}

\bibitem{KT1} R. KOMOROWSKI and N. TOMCZAK-JAEGERMANN, {\em Banach spaces
without local unconditional structure}, Israel Journal of Mathematics {\bf 89}
(1995), 205-226.

\bibitem{KT2} R. KOMOROWSKI and N. TOMCZAK-JAEGERMANN, {\em Erratum to:
"Banach spaces without local unconditional structure"},
 Israel Journal of Mathematics {\bf 105}
(1998), 85-92.

\bibitem{lt2}J. LINDENSTRAUSS and L. TZAFRIRI: Classical Banach Spaces II,
\emph{Springer,
Berlin, Heidelberg 1979}.

















\bibitem{T} N. TOMCZAK-JAEGERMANN, {\em Banach spaces with many
  isomorphisms}, Semin\'ario Brasileiro de An\'alise {\bf 48}, Petropolis,
  Rio de Janeiro,  November 1998, 189-210.

\end{thebibliography}
\end{document}